\documentclass{amsart}
 \usepackage{amscd}
                     \usepackage{amssymb}
 \newcommand{\be}{\begin{equation}}
       \newcommand{\ee}{\end{equation}}
       \newcommand{\ba}{\begin{eqnarray}}
        \newcommand{\ea}{\end{eqnarray}}
 \newcommand{\ban}{\begin{eqnarray*}}
 \newcommand{\ean}{\end{eqnarray*}}

\def\XXint#1#2#3{{\setbox0=\hbox{$#1{#2#3}{\int}$}
     \vcenter{\hbox{$#2#3$}}\kern-.5\wd0}}

 \newcommand{\lp}{\langle}
 \newcommand{\rp}{\rangle}
 \newcommand{\ra}{\rightarrow}

  \newcommand{\Pf}{\noindent {\bf Proof:} }
  \newcommand{\Rk}{\noindent {\bf Remark} }
 
 \newcommand{\sect}[1]{\section{#1} \setcounter{equation}{0}}

 \newcommand{\vol}{\mathrm{Vol}}
 
 \newcommand{\Ric}{\mathrm{Ric}}
 
 \newcommand{\Hess}{\mathrm{Hess}}
 \newcommand{\divg}{\mathrm{div}}

  \newcommand{\tr}{\mbox{tr}}

 \newtheorem{theo}{Theorem}[section]

\begin{document}
 \newtheorem{defn}[theo]{Definition}
 \newtheorem{ques}[theo]{Question}
 \newtheorem{lem}[theo]{Lemma}
 \newtheorem{lemma}[theo]{Lemma}
 \newtheorem{prop}[theo]{Proposition}
 \newtheorem{coro}[theo]{Corollary}
 \newtheorem{ex}[theo]{Example}
 \newtheorem{note}[theo]{Note}
 \newtheorem{conj}[theo]{Conjecture}
 \title{Monotonicity Formulas for Bakry-Emery Ricci Curvature}
\author{Bingyu Song}
 \thanks{Partially  supported by  NSFC Grant No. 11171259}
 \email{bysong@mail.ccnu.edu.cn}
 \address{ School of Mathematics and Statistics,  Central China Normal University, Wuhan 430079, China}
\author{Guofang Wei}
 \thanks{Partially supported by NSF Grant \# DMS-1105536.}
 \email{wei@math.ucsb.edu}
 \address{Department of Mathematics,  UCSB,  Santa Barbara, CA 93106}
 \author{ Guoqiang Wu}
\email{cumtwgq@mail.ustc.edu.cn}
\address{Department of Mathematics, University of Science and Technology of China}
 \date{}
 \maketitle

 \begin{abstract}Motivated and inspired by the recent work of Colding \cite{Colding} and Colding-Minicozzi \cite{Colding-Minicozzi} we derive several families of monotonicity formulas for manifolds with nonnegative Bakry-Emery Ricci curvature, extending the formulas in \cite{Colding,Colding-Minicozzi}.
  \end{abstract}

 \sect{Introduction}
The Bakry-Emery Ricci tensor is a Ricci tensor for smooth metric
measure spaces, which are Riemannian manifolds with  measures
conformal to the Riemannian measures. Formally a smooth metric
measure space is a triple $(M^n,g, e^{-f} dvol_g)$, where $M$ is a
complete $n$-dimensional Riemannian manifold with metric $g$, $f$ is
a smooth real valued function on $M$, and $dvol_g$ is the Riemannian
volume density on $M$. These spaces occur naturally as smooth collapsed limits
of manifolds  under the measured Gromov-Hausdorff convergence.

The $N$-Bakry-Emery Ricci tensor is
\be
\Ric_f^N = \Ric +\Hess f  -\frac 1N df \otimes df  \ \ \  \mbox{for} \   0\le N \le \infty.   \label{Ric-N-f}
\ee
The purported  dimension of the  space is related to $N$, i.e. it is $n+N$.  When $N=0$, we assume $f$ is constant and  $\Ric_f^N = \Ric$, the usual Ricci curvature.
When $N$ is infinite,  we denote
$\Ric_f = \Ric_f^\infty = \Ric +\Hess f$. Note that if $N_1 \geq N_2$ then $\Ric_f^{N_1}
\geq \Ric_f^{N_2}$ so $\Ric_f^N \ge \lambda g$ implies $\Ric_f
\ge \lambda g.$

The Einstein equation $\Ric_f = \lambda g$ ($\lambda$ a constant)  is exactly
the gradient Ricci soliton equation, which plays an important role in the theory
of Ricci  flow. On the other hand,  the equation $\Ric_f^N = \lambda g$, for $N$ a positive integer, corresponds to
warped product Einstein metric on $M \times_{e^{-\frac fN}} F^N$, where $F^N$ is  some $N$ dimensional Einstein manifold, see \cite{Case-Shu-Wei}.

Recently Colding \cite{Colding} and Colding-Minicozzi \cite{Colding-Minicozzi}  introduced some new monotonicity formulas associated to positive Green's function of the Laplacian. These formulas are very useful and related to other known monotonicity formulas, see \cite{Colding-Minicozzi-survey}. In particular, using one of the monotonicty formula Colding-Minicozzi showed that for any Ricci-flat manifold with Euclidean volume growth, tangent cones at
infinity are unique as long as one tangent cone has a smooth cross-section \cite{Colding-Minicozzi-cone}. Our paper is motivated and inspired by these work of Colding and Colding-Minicozzi.

With respect to the measure $e^{-f} dvol$ the natural self-adjoint  $f$-Laplacian is  $\Delta_f = \Delta - \nabla f \cdot  \nabla$. Consider  the positive Green's function $G(x_0, \cdot)$ of the $f$-Laplacian of $(M^n, g, e^{-f} dvol)$ (see Definition~\ref{Green}).
For any real number $k > 2$,
let $b=G^{\frac{1}{2-k}}$.
For $\beta,l, p \in \mathbb R$, when $b$ is proper, we consider
$$A^{\beta}_f(r)=r^{1-l}\int_{b=r}|\nabla b|^{\beta+1}e^{-f},$$
 \[ V^{\beta,p}_f(r)=r^{p-l}\int_{b\leq r}\frac{|\nabla b|^{2+\beta}}{b^{p}}e^{-f}.\]
While $A^{\beta}_f(r)$ is well defined for all $r>0$, $V^{\beta,p}_f(r)$ is only well defined when \be
C(n,k,p) = (n-2)(k-p) -\beta(k-n) > 0. \label{C-n-k-p}
\ee  See the proof of Lemma~\ref{A-V-0}  for detail.
When $k=l=n, \beta =2, p=0$, these reduce to $A(r), V(r)$ in \cite{Colding}. When $k=l=n, p=2$, these are $A_\beta, V_\beta$ in \cite{Colding-Minicozzi}.

First we obtain the following gradient estimate for $b$.
\begin{prop} \label{nabla-b-p}
If a smooth metric measure space $(M^n, g, e^{-f} dvol)$ ($n \ge 3$) has $\Ric_f^N \ge 0$, then for $k = n+N$, there exists $r_0 >0$, such that on $M\setminus B(x_0, r_0)$,
\be
|\nabla b(y) | \le  C(n, N, r_0).  \label{nabla-b}
\ee
\end{prop}
\Rk In  \cite[Theorem 3.1]{Colding} Colding obtained the sharp estimate that if $\Ric_{M^n} \ge 0$ ($n \ge 3$), then $|\nabla b| \le 1$ for $f=0, k=n$ in above. From (\ref{b-0}) this can not be true when $k>n$ as $|\nabla b(y)| \ra \infty$ as $y\ra x_0$.  For $\Ric_f \ge 0$, $|\nabla b(y)|$ may not be bounded as $y \ra \infty$, see Example~\ref{Bryant-soliton}.

We prove many families of monotonicity formulas, which, besides  recovering the ones in \cite{Colding, Colding-Minicozzi} when $k=l=n$ and $f$ is constant,  give some new ones even in this case. For example, when $N$ is finite, we have

\begin{theo} \label{A-V-mono}
If $M^n (n\ge 3)$ has $\Ric^N_f \ge 0$, then, for $k \ge n+N$, $k\leq l\leq 2k-2$, $\alpha = 3k-p-l-2$ and $C(n,k,p) >0$,
\ba
\lefteqn{(A^{\beta}_f - \alpha V^{\beta,p}_f)'(r)} \nonumber \\
& &  \ge r^{p-1-l} \int_{b\leq r} \frac{\beta |\nabla b|^{\beta -2}}{4b^{p} }
\bigg\{ \left|\Hess\,b^2-\frac{\Delta b^2}{n}g\right|^2  +4(\beta-2)b^2|\nabla|\nabla b||^2\bigg\}e^{-f}. \nonumber \\
& &
\ea
Hence if in addition $\beta \ge 2$, then $A^{\beta}_f - \alpha V^{\beta,p}_f$ is nondecreasing in $r$.
\end{theo}

See discussion in Section~5 for full generality. Note that when $N=0$ (i.e. $f$ is constant and $\Ric \ge 0$), and $p=0, \beta \ge 2$, $k=n$ we get monotonicity for all $n \le l \le 2n-2$;  in the case  when $\beta = 2, l =n$, this is the first monotonicity formula in \cite{Colding}.

\begin{theo}  \label{A-mono}
If $M^n (n\ge 3)$ has $\Ric^N_f \ge 0$, then for $\beta \ge 2$, $k=l= n+N$, $(A^\beta_f)'(r) \le 0$ and $(V_f^{\beta, p})'(r) \le 0$ for $p<n+N-\frac{\beta N}{n-2}$. In fact \be
(A^\beta_f)'(r) \le - \frac{\beta}{4} r^{k-3} \int_{ b\ge r} b^{2-2k}|\nabla b|^{\beta-2} \left|\Hess\,b^2-\frac{\Delta b^2}{n}g\right|^2 e^{-f}.  \label{A'k=l}
\ee
\end{theo}
Again this reduces to a formula in \cite{Colding} when $k= l =n$ and $f$ is constant, which is used in \cite{Colding-Minicozzi-cone} to show that  for any Ricci-flat manifold with Euclidean volume growth, tangent cones at
infinity are unique as long as one tangent cone has a smooth cross-section.

Without any assumption on $f$ we also establish monotonicity formulas when $N$ is infinite.  For example,
\begin{theo} \label{A-V-mono-infty}
If $M^n (n\ge 4)$ has $\Ric_f \ge 0$,  then for $\beta =2, p=0, k \ge 12, l=\frac 32 k -1$,  we have $(A^{\beta}_f - (\frac 32 k-1) V^{\beta,p}_f)'(r) \ge 0$;
for $\beta =2, k \ge 12, l = \frac 32(k-1)$, $r_2 \ge r_1 >0$, we have $(A_f^\beta)' (r_2) \ge (A_f^\beta)' (r_1)$.
\end{theo}
Here $l>k$. See Theorem~\ref{A-V-mono-infty-gen} for more general statement.

When $k=l$, similar monotonicity is not true any more. In fact for $k=l=n$ we show while for $n=3$ several monotonity still holds for the Bryant solitons for $r$ large, it does not hold when $n \ge 5$.  In fact it is monotone in the opposite direction for $r$ large,  see Example~\ref{Bryant-soliton}.   With some condition on $f$ we still get several monotonicity, see e.g. Corollary~\ref{mono-k=l=n}.

As in \cite{Colding-Minicozzi}, one can use the term $\left|\Hess\,b^2-\frac{\Delta b^2}{n}g\right|^2$ to allow smaller $\beta$ (one only needs $\beta \ge 1-\frac{1}{n-1}$ instead of $\beta\ge 2$, see the last part of Section 5).

The paper is organized as follows. In the next section, we discuss the existence and the basic properties of Green's function for $f$-Laplacian and prove Proposition~\ref{nabla-b-p}. In Section 3, using the following Bochner formula (see e.g. \cite{Wei-Wylie2007})
\be  \label{Bochner-f}
\frac 12 \Delta_f
|\nabla u|^2 = |\Hess \, u|^2 + \lp \nabla u, \nabla (\Delta_f u) \rp
+ \Ric_f (\nabla u, \nabla u).
\ee
we compute the f-laplacian of $(b^{2q} |\nabla b|^{\beta})$, a key formula needed for deriving the monotonicity formulas in Section 4. In Section 5, 6 we apply these formulas to the case when $\Ric_f^N \ge 0$ for $N$ finite and infinite respectively.

As with many other monotonicity formulas, we expect our formulas will have nice applications, especially for quasi-Einstein manifolds and steady gradient Ricci solitons.

{\bf Acknowledgments}.   The authors would like to thank Toby Colding for his interest and encouragement, Toby Colding and Bill Minicozi for answering our questions on their work. This work was done while the first and third authors were visiting UCSB. They would like to thank UCSB for  hospitality during their stay.

\sect{Green's function}

\begin{defn}  \label{Green}
Given a smooth metric measure space $(M^n, g, e^{-f} dvol)$, $x_0\in M$, $G = G(x_0, \cdot)$ is the Green's function of  the $f$-Laplacian (with pole at $x_0$)  if
\[
\Delta_f G = - \delta_{x_0}.\]
\end{defn}

In \cite{Malgrange1955, Li-Tam1987} it is shown that on any complete Riemannian
manifold there exists a symmetric Green's function of the Laplacian.  Same proof carries over for $f$-Laplacian.

$M$ is called $f$-nonparabolic if it has positive Green's function. When $f=0$ and $\Ric_M \ge 0$, the existence of positive Green's function is well understood \cite{Varopoulos1981}, see also \cite{Li-Yau1986, Sung1998}. Namely $M$ is nonparabolic if and only if $\int_1^\infty  \frac{r}{\vol B(x_0, r)} dr < \infty$, and $G \ra 0$ at infinity. Same result also holds when $\Ric_f^N \ge 0$ ($N$ finite) or $\Ric_f \ge 0$ and $f$ is bounded. This and other existence results will be studied in \cite{Dai-Sung-Wang-Wei}. In there we observe that all nontrivial steady Ricci solitons are $f$-nonparabolic.

In this paper we assume $M$ is $f$-nonparabolic and $G \ra 0$ at infinity so $b$ is proper. Also, near the pole (after normalization),
for $n \ge 3$, by \cite{Gilbarg-Serrin1955}, we have \be
G(y) = d^{2-n} (x_0,y) (1+o(1)), \ \ |\nabla G(y)|  = (n-2) d^{1-n}(x_0,y)(1+o(1)),  \label{G-0}  \ee
where $o(1)$ is a function with $o(y) \ra 0$ as $y \ra x_0$.

Recall for $k >2$,  $b=G^{\frac{1}{2-k}}$. Then
\be b(y) = d^{\frac{n-2}{k-2}}(x_0,y) (1+o(1)), \ \ |\nabla b(y)| = \frac{n-2}{k-2} d^{\frac{n-k}{k-2}}(x_0,y) (1+o(1)). \label{b-0}
\ee

  To prove Proposition~\ref{nabla-b-p}, we need the following two propositions for $\Ric_f^N \ge 0$. Similar to the case when $\Ric \ge 0$,  one has the following Laplacian comparison and gradient estimate for $\Ric_f^N \ge 0$  \cite{Qian1995, Qian1997}, see also \cite{Wei-Wylie2007, Li2005}.
Let $r(x) = d(x, x_0)$ be the distance function, then
\be
\Delta_f r \le \Delta_{\mathbb R^{n+N}} r = \frac{n+N}{r}. \label{lap-com}
\ee
If $u$ is a $f$-harmonic function on $B(x, R)$, then on $B(x, R/2)$,
\be
|\nabla \log u| \le \frac{C(n+N)}{R}.  \label{gradient}
\ee

\noindent {\bf Proof of Proposition~\ref{nabla-b-p}}.
 For any $y \in M\setminus \{x_0\}$, $G$ is a smooth harmonic function on $B(y, r)$ with $r = r(y)=d(y, x_0)$. By (\ref{gradient}),
 \[
 |\nabla \log G | \le \frac{C(k)}{r} \ \ \ \mbox{on} \ \ B(y, r/2).
 \]
 Now $|\nabla \log b| = \frac{1}{k-2}  |\nabla \log G |$. Therefore \[
 |\nabla b| (y) =  |\nabla \log b| \cdot b \le  \frac{C(k)}{k-2} \cdot  \frac br.
 \]
 Hence (\ref{nabla-b}) follows if we show  $b  \le   C_1(n,N,r_0) r$, i.e. $G \ge C_2(n,N,r_0) r^{2-k}$ for some $r_0 >0$ on $M\setminus B(x_0, r_0)$.

For any $\epsilon >0$,  we have on $M\setminus \{x_0\}$, \[
 \Delta_f (G - \epsilon\,  r^{2-k}) = -\epsilon \, \Delta_f ( r^{2-k} ) \le -\epsilon \,  \Delta_{\mathbb R^k} (r^{2-k} ) =0,\]
 where the inequality follows from the Laplacian comparison (\ref{lap-com}).
 Since $\lim_{r \ra \infty} (G - \epsilon\,  r^{2-k})  = 0$. Also by (\ref{G-0}),  there exists $r_0 >0$ small such that $G(y) \ge \frac 12 r_0^{2-n}$ on $\partial B (x_0,r_0)$.   Take $\epsilon = \frac 12 r_0^{k-n}$, we have  $(G - \frac 12 r_0^{k-n} r^{2-k}) \ge 0$ on $\partial B (x_0,r_0)$.    Therefore
by the maximum principle we have $G(y) \ge   \frac 12 r_0^{k-n} r^{2-k}$ when $r(y) \ge r_0$. Namely $b  \le  \left(2 r_0^{n-k}\right)^{\frac{1}{k-2}} r$.
\qed

\section{The $f$-Laplacian of $b$ and $|\nabla b|$}

\begin{lem} For any real number $\beta$,
\be
\Delta_f b=\frac{k-1}{b} |\nabla b|^2,   \label{Delta-f-b}
\ee
\be
\Delta_f b^\beta=\beta(\beta + k-2)b^{\beta-2}  |\nabla b|^2.  \label{Delta-f-b-alpha}
\ee
In particular,
\be
\Delta_f b^2= 2k  |\nabla b|^2.  \label{Delta-f-b-2}
\ee
\end{lem}
\Pf For any positive function $v$, we have
\be
\Delta_f v^{\beta}
=\beta v^{\beta -1} \left[(\beta-1) \frac{|\nabla v|^2}{v} +\Delta_f v \right]  \label{v-alpha}
\ee
Since $\Delta_f b^{2-k} = 0$, this gives $(2-k-1)\frac{|\nabla b|^2}{b} +\Delta_f b=0$, namely (\ref{Delta-f-b}).
Combining (\ref{v-alpha}) and (\ref{Delta-f-b}) gives (\ref{Delta-f-b-alpha}).
\qed

The following important formulas holds for any positive $f$-harmonic function $G$, not just Green's function.
\begin{prop}
\ba
\Delta_f |\nabla b|^\beta &=&\frac{\beta}{4b^2}|\nabla b|^{\beta-2}
\left\{ \left|\Hess\, b^2\right|^2
+\Ric_f (\nabla b^2, \nabla b^2)+2(k-2)  \lp \nabla b^2, \nabla |\nabla b|^2 \rp \right. \nonumber\\
& &+4(\beta -2)b^2 |\nabla|\nabla b||^2 -4k|\nabla b|^4 \Big\}.  \label{Delta-f-grad-b-alpha}
\ea
\ba
\lefteqn{ \Delta_f (b^{2q} |\nabla b|^{\beta})   }   \nonumber  \\
&  = &  \frac{\beta}{4} b^{2q-2} |\nabla b|^{\beta -2}
\bigg\{ \left|\Hess\, b^2\right|^2 + \Ric_f (\nabla b^2, \nabla b^2)
 +2(k-2+2q)  \lp \nabla b^2, \nabla |\nabla b|^2 \rp   \nonumber  \\
 & &  +4(\beta -2) b^2 |\nabla |\nabla b||^2 +  \left[\frac{8q}{\beta}(k-2+2q)-4k\right]|\nabla b|^4
\bigg\} \label{Delta-f-b-beta}
\ea
\end{prop}

\Pf  Since the formulas are in terms of the function $b^2$, we first compute $ \Delta_f |\nabla b^2|^2$. Applying  the Bochner formula (\ref{Bochner-f}) to $b^2$ and using (\ref{Delta-f-b-2}),  we have
\ba
\frac 12 \Delta_f |\nabla b^2|^2  &  =  &  |\Hess\, b^2|^2+\Ric_f (\nabla b^2, \nabla b^2)  +  \lp \nabla b^2, \nabla (\Delta_f b^2) \rp   \nonumber \\
&  =  &  |\Hess\, b^2|^2 +\Ric_f (\nabla b^2, \nabla b^2) + 2k \lp \nabla b^2, \nabla |\nabla b|^2 \rp .   \label{Delta-f-grad-b2}
\ea
Now we  compute $\Delta_f |\nabla b|^2$.  Since $|\nabla b^2|^2 = 4b^2 |\nabla b|^2$,
\ban
\frac 12 \Delta_f |\nabla b^2|^2  &   =  & 2 \left( |\nabla b|^2 \Delta_f b^2
+  b^2 \Delta_f |\nabla b|^2 + 2 \lp \nabla b^2, \nabla |\nabla b|^2 \rp  \right) \\
& = & 4k |\nabla b|^4 + 2b^2 \Delta_f |\nabla b|^2 + 4  \lp \nabla b^2, \nabla (|\nabla b|^2) \rp.
\ean
Combine this with (\ref{Delta-f-grad-b2}),  we have
\ba
2b^2 \Delta_f |\nabla b|^2  &  =  &  |\Hess\, b^2|^2 +\Ric_f (\nabla b^2, \nabla b^2)
+ (2k-4) \lp \nabla b^2, \nabla |\nabla b|^2 \rp   \nonumber  \\
& &-4k|\nabla b|^4  .
\ea
Then,
\ban
\Delta_f |\nabla b|^{\beta} &=& \Delta_f (|\nabla b|^2)^{\frac{\beta}{2}} \\
&=& \frac{\beta}{2}|\nabla b|^{\beta -2}  \Delta_f |\nabla b|^2
+\beta (\beta -2)|\nabla b|^{\beta-2} |\nabla |\nabla b||^2 \\
&=&\frac{\beta}{4b^2} |\nabla b|^{\beta -2}
\left\{ |\Hess\, b^2|^2 + \Ric_f (\nabla b^2, \nabla b^2)+ (2k-4)  \lp \nabla b^2, \nabla |\nabla b|^2 \rp  \right.   \nonumber  \\
& & \left. +4(\beta -2) b^2 |\nabla |\nabla b||^2  -4k|\nabla b|^4 \right\},
\ean
which is  (\ref{Delta-f-grad-b-alpha}).

For the second one,  by the product formula for Laplacian and using (\ref{Delta-f-b-alpha}), we get
\ba
\Delta_f (b^{2q} |\nabla b|^{\beta}) &=& b^{2q} \Delta_f |\nabla b|^{\beta}
+|\nabla b|^{\beta}  \Delta_f (b^{2q}) +2 \lp \nabla |\nabla b|^{\beta}, \nabla b^{2q} \rp  \nonumber \\
&=&  b^{2q} \Delta_f |\nabla b|^{\beta}
+2q(2q+k-2)b^{2q-2} |\nabla b|^{2+\beta}  \nonumber \\
& & + \beta q |\nabla b|^{\beta -2} b^{2q-2}  \lp \nabla b^2, \nabla |\nabla b|^2 \rp. \nonumber
\ea
Plug in (\ref{Delta-f-grad-b-alpha}) we obtain (\ref{Delta-f-b-beta}).   \qed

\sect{Monotonicity Formulas}

Recall for $l,\beta, p \in \mathbb R$,
$$A^{\beta}_f(r)=r^{1-l}\int_{b=r}|\nabla b|^{\beta+1}e^{-f},$$
 \[ V^{\beta,p}_f(r)=r^{p-l}\int_{b\leq r}\frac{|\nabla b|^{2+\beta}}{b^{p}}e^{-f}.\]

As $r \ra 0$, we have the following information.
\begin{lem}  \label{A-V-0}
Let $M^n$ be a smooth manifold with $n \ge 3$. Denote $C(n,k,l) = (k-l)(n-2)+(n-k)\beta$. Then
\ba
\lim_{r \ra 0} A_f^\beta (r) & = & \left\{ \begin{array}{ll} 0 &  \mbox{if} \ C(n,k,l) >0 \\
\left( \frac{n-2}{k-2} \right)^{1+\beta} Vol(\partial B_1(0))e^{-f(x_0)} & \mbox{if} \ C(n,k,l) = 0 \\
\infty & \mbox{if} \ C(n,k,l) < 0 \end{array}, \right. \label{A-0}\\
\lim_{r \ra 0} V_f^{\beta, p} (r) & = & \left\{ \begin{array}{ll} 0 &  \mbox{if} \ C(n,k,l) >0 \\
\left( \frac{n-2}{k-2} \right)^{1+\beta} \frac{n-2}{C(n,k,p)}\vol(\partial B_1(0)) e^{-f(x_0)} & \mbox{if} \ C(n,k,l) = 0 \\
\infty  & \mbox{if} \ C(n,k,l) < 0 \end{array}, \right.  \nonumber  \\
\label{V-0}
\ea
where $\vol (\partial B_1(0))$ is  the volume of the unit sphere in $\mathbb R^n$.
\end{lem}
\Pf From (\ref{b-0}),
\ba
A_f^\beta(r)&=&r^{1-l}\int_{b=r}|\nabla b|^{1+\beta}e^{-f}\nonumber\\
            &=&\left( \frac{n-2}{k-2} \right)^{1+\beta} \left(r^{1-l+\frac{n-k}{n-2}(1+\beta)+\frac{k-2}{n-2}(n-1)}\right)(1+o(1))e^{-f(x_0)} \vol(\partial B_1(0)), \nonumber
\ea
where $o(1) \ra 0$ as $r \ra 0$.
Note that $1-l+\frac{n-k}{n-2}(1+\beta)+\frac{k-2}{n-2}(n-1) = (k-l)(n-2)+(n-k)\beta$. This gives (\ref{A-0}).

Similarly,
\ba
V_f^{\beta,p}(r)&=&r^{p-l}\int_0^r\int_{b=s}\frac{|\nabla b|^{1+\beta}}{b^p}e^{-f}\nonumber\\
&=& \left( \frac{n-2}{k-2} \right)^{1+\beta} r^{p-l}\int_0^r s^{\frac{n-k}{n-2}(1+\beta)-p+\frac{k-2}{n-2}(n-1)}ds(1+o(1))\cdot Vol(\partial B_1(0))e^{-f(x_0)}\nonumber
\ea
The integral exists if the constant in (\ref{C-n-k-p}), $C(n,k,p)>0$, and
\ba
V_f^{\beta,p}(r)=\frac{n-2}{C(n,k,p)}\left( \frac{n-2}{k-2} \right)^{1+\beta}r^{\frac{C(n,k,l)}{n-2}}(1+o(1))Vol(\partial B_1(0))e^{-f(x_0)}\nonumber
\ea
This gives (\ref{V-0}).
\qed

Since
\[
V^{\beta,p}_f (r)=r^{p-l}\int_0^r \int_{b=s} \frac{|\nabla b|^{1+\beta}}{b^{p}}e^{-f}.
\]
we have
\ba
(V_f^{\beta,p})' (r) &=& (p-l) r^{p-l-1} \int_0^r \int_{b=s} \frac{|\nabla b|^{1+\beta}}{b^{p}}e^{-f} +r^{p-l} \int_{b=r}  \frac{|\nabla b|^{1+\beta}}{b^{p}}e^{-f}  \nonumber \\
&=& \frac{p-l}{r} V^{\beta,p}_f (r) +\frac 1r A^{\beta}_f (r).  \label{V'}
\ea

To find the derivative of $A_f (r)$,  we use the following formula.
\begin{lem}
For a smooth function $u: M\setminus x \ra \mathbb{R} $, let
\[
I_u (r)= \int_{b=r}u |\nabla b|\, e^{-f}.
\]
Then for any $0 < r_0 \le r$,
\ba
I'_u (r)&=& \frac{k-1}{r} I_u (r) + \int_{b=r}\lp \nabla u, \nu \rp e^{-f} \label{Iu'} \\
&=& \frac{k-1}{r} I_u (r) + \int_{r_0 \le b\leq r} \left( \Delta_f u \right) e^{-f}+  \int_{b=r_0}\lp \nabla u, \nu \rp e^{-f},
  \label{Iu}
\ea
where $\nu =  \frac{\nabla b}{|\nabla b|}$ is the unit normal direction.
\end{lem}
This formula can be derived using the diffeomorphisms  generated by $\frac{\nabla b}{|\nabla b|^2}$, see \cite[Appendix]{Colding-Minicozzi1997}. For completeness, we give a simple proof using just the divergent theorem and the co-area formula.

\Pf   Note that
\ba
I_u (r) - I_u (r_0)  &=&\int_{b=r}u|\nabla b| \, e^{-f} - \int_{b=r_0}u|\nabla b| \, e^{-f}
\nonumber\\
&=&\int_{r_0 \le b\leq r} \divg(u e^{-f}\nabla b)\nonumber\\
&=&\int_{r_0}^r\int_{b=s}\frac{\divg (u e^{-f}\nabla b)}{|\nabla b|}.\nonumber
\ea
Take derivative of this equation both sides with respect to $r$ and using (\ref{Delta-f-b}) gives
\ban
I_u'&=&\int_{b=r}\frac{\divg(u e^{-f}\nabla b)}{|\nabla b|}\nonumber\\
&=&\int_{b=r}\frac{\lp \nabla u, \nabla b \rp  e^{-f}+u e^{-f}\Delta_f b}{|\nabla b|}\nonumber\\
&=&\int_{b=r}\lp \nabla u,  \frac{\nabla b}{|\nabla b|} \rp e^{-f}+\int_{b=r}u e^{-f}\frac{k-1}{b}|\nabla b|\nonumber \\
& =& \int_{r_0 \le b\leq r} \left( \Delta_f u \right) e^{-f} + \int_{b=r_0}\lp \nabla u,  \frac{\nabla b}{|\nabla b|} \rp e^{-f}+\frac{k-1}{r}\int_{b=r}u |\nabla b|\, e^{-f}.
\ean
\qed

In order to match the derivative of $A_f^\beta$ with $V_f^{\beta, p}$, we write $A_f^{\beta}(r) = r^{p-l-1} \displaystyle\int_{b=r}b^{2-p}|\nabla b|^{\beta+1}e^{-f}$. Applying (\ref{Iu}) to $u = b^{2-p}|\nabla b|^{\beta}$,  if $C(n,k,p) > 0$, the integral on $b=r_0$ goes to zero as $r_0 \ra 0$, and we have
\begin{coro}
\be
(A_f^{\beta})'(r) =\frac{k-l-2+p}{r} A_f^{\beta}(r) +r^{p-1-l} \int_{b \leq r} \left(\Delta_f(b^{2-p} |\nabla b|^\beta) \right) e^{-f}.  \label{A'}
\ee
\end{coro}

As in \cite{Colding} we derive a formula which will give  the first monotonicity.
\begin{theo} \label{A-V'} When $k>2$,  $C(n,k,p) > 0$, for any $\alpha  \in \mathbb R$,
\ba
\lefteqn{(A^{\beta}_f -\alpha V^{\beta,p}_f)'(r) = }   \nonumber\\
& & r^{p-1-l} \int_{b\leq r} \frac{\beta |\nabla b|^{\beta -2}}{4b^{p} }
\bigg\{ \left|\Hess\, b^2\right|^2 + \Ric_f (\nabla b^2, \nabla b^2) +4(\beta-2)b^2|\nabla|\nabla b||^2\bigg\}e^{-f}\nonumber\\
& & + \frac{1}{r} \left( \lambda_1 A_f^\beta (r) + \lambda_2  V_f^{\beta,p}(r) \right),
\label{A-V}
\ea
where \ba
\lambda _1 & = & 3k-p-l-2-\alpha, \\
\lambda_2  & =&  (p+2-2k)(k-p)-\beta k -\alpha(p-l).
\ea
\end{theo}

\Pf  From (\ref{A'}) we would like to compute $\displaystyle\int_{b \leq r} \left(\Delta_f(b^{2-p} |\nabla b|^\beta) \right) e^{-f}.$
By (\ref{Delta-f-b-beta}),
\ba
\Delta_f (b^{2-p} |\nabla b|^{\beta}) &=&\frac{\beta |\nabla b|^{\beta -2}}{4b^{p}}
\Bigg\{ \left|\Hess\, b^2\right|^2 + \Ric_f (\nabla b^2, \nabla b^2)
+2(k-p)  \lp \nabla b^2, \nabla |\nabla b|^2 \rp \nonumber  \\
& & +4(\beta -2) b^2 |\nabla |\nabla b||^2 +  \left[\frac{8-4p}{\beta}(k-p)-4k\right]|\nabla b|^4
\Bigg\}.  \label{Delta-f-b-p}
\ea
To compute the third term, by Stokes' theorem and (\ref{Delta-f-b-alpha}),
\ba
\lefteqn{ \int_{b\leq r}  b^{-p} \frac{\beta}{2} |\nabla b|^{\beta -2} \lp \nabla b^2, \nabla|\nabla b|^2 \rp e^{-f} } \nonumber \\
& = & \int_{b\leq r} \frac{2}{2-p} \lp  \nabla b^{2-p}, \nabla \left(|\nabla b|^2 \right)^{\beta/2} \rp e^{-f} \nonumber\\
&=&\int_{b=r}b^{-p}\lp \nabla b^2,\frac{\nabla b}{|\nabla b|} \rp |\nabla b|^\beta e^{-f}- \frac{2}{2-p} \int_{b\leq r} \Delta_f (b^{2-p}) |\nabla b|^{\beta}e^{-f}  \nonumber\\
&=&2r^{1-p}\int_{b=r}|\nabla b|^{\beta+1}e^{-f}-2(k-p)\int_{b\leq r}b^{-p} |\nabla b|^{\beta+2} e^{-f} \nonumber\\
&=&2r^{l-p} A_f^\beta +2(p-k)r^{l-p}V_f^{\beta,p}. \label{nabla-nabla}
\ea
In the above proof we assume $p \not= 2$. By taking limit of (\ref{nabla-nabla}) as $p \ra 2$, we see (\ref{nabla-nabla}) also holds for $p=2$. In the second equality, we use $C(n,k,p) > 0$ so the integral for $r \ra 0$ is zero.

Combining (\ref{nabla-nabla}),(\ref{Delta-f-b-p}), (\ref{A'}) and (\ref{V'}), we have
\ba
\lefteqn{(A^{\beta}_f -\alpha V^{\beta,p}_f)'(r) = } \nonumber\\
& & r^{p-1-l} \int_{b\leq r} \frac{\beta |\nabla b|^{\beta -2}}{4b^{p} }
\bigg\{ \left|\Hess\, b^2\right|^2 + \Ric_f (\nabla b^2, \nabla b^2) +4(\beta-2)b^2|\nabla|\nabla b||^2\bigg\}e^{-f}\nonumber\\
& &+ \frac{3k-p-l-2-\alpha}{r}A_f^\beta (r) + \frac{1}{r}\left[(p+2-2k)(k-p)-\beta k -\alpha(p-l)\right]V_f^{\beta,p}(r).\nonumber
\ea
\qed

Letting $k=l, \alpha = 2k-p-2$ in (\ref{A-V}) gives
\begin{coro}
 When $k>2$,  $C(n,k,p) > 0$, and $k=l$, we have
\ba
\lefteqn{(A^{\beta}_f -(2k-p-2) V^{\beta,p}_f)'(r) = r^{p-1-l} \int_{b\leq r} \frac{\beta |\nabla b|^{\beta -2}}{4b^{p} }
\bigg\{ \left|\Hess\, b^2\right|^2  }   \nonumber\\
& & + \Ric_f (\nabla b^2, \nabla b^2) +4(\beta-2)b^2|\nabla|\nabla b||^2 - 4k |\nabla b|^4 \bigg\}e^{-f}.  \label{A-V'-k=l}
\ea
\end{coro}

Following is a formula which will give second monotonicity formula.
\begin{theo} \label{rA'}  For $c, d \in \mathbb R$, let $g(r) = r^c \left(r^d A_f^\beta (r) \right)'$. Then for $0< r_1 < r_2$,
\ba
\lefteqn{ g(r_2)- g(r_1)} \nonumber \\
 & =&\int_{r_1\leq b\leq r_2}\frac{\beta}{4}b^{c+d-l-1}|\nabla b|^{\beta-2}\bigg\{|\Hess\, b^2|^2+\Ric_f(\nabla b^2,\nabla b^2) + \lambda_3|\nabla b|^4
\nonumber \\
&  + & 4(\beta-2)b^2|\nabla|\nabla b||^2 \bigg\}e^{-f} + \lambda_4 \int_{r_1\leq b\leq r_2}b^{c+d-l}\lp \nabla b, \nabla|\nabla b|^\beta\rp e^{-f}, \label{A'-mono}
\ea
where
\ba
\lambda_3 & = & \frac{4}{\beta}(k+d-l+c-1)(k+d-l) -4k, \\
\lambda_4 & = & 3k-2l-3+c+2d.
\ea
\end{theo}
\Pf From (\ref{Iu'}) with $u = |\nabla b|^\beta$, we have
\ba
(A_f^\beta)'(r) =r^{1-l}\int_{b=r}\lp \nabla|\nabla b|^\beta,\nu \rp e^{-f} +
(k-l)r^{-l}\int_{b=r}|\nabla b|^{1+\beta}e^{-f}.\nonumber
\ea
Hence
\ba
g(r) =\int_{b=r}b^{c+d+1-l}\lp \nabla|\nabla b|^\beta,\nu \rp e^{-f} +
(k+d-l)\int_{b=r} b^{c+d-l} |\nabla b|^{1+\beta}e^{-f}, \nonumber
\ea
and
\ba
g(r_2)-g(r_1)
&=&\int_{r_1\leq b\leq r_2}\divg(b^{c+d+1-l} e^{-f}\nabla|\nabla b|^\beta) \nonumber \\
& & +(k+d-l)\int_{r_1\leq b\leq r_2}\divg(b^{c+d-l}|\nabla b|^\beta e^{-f} \nabla b). \nonumber
\ea
Since
\ban
\lefteqn{ \divg(b^{c+d+1-l}e^{-f} \nabla|\nabla b|^\beta)} \\
& = & (c+d+1-l)b^{c+d-l}\lp \nabla b, \nabla|\nabla b|^\beta\rp e^{-f}+b^{c+d+1-l}(\Delta_f|\nabla b|^\beta) e^{-f}  \ean
 and
 \ban
 \divg(b^{c+d-l}|\nabla b|^\beta e^{-f} \nabla b)& = & (c+d-l)b^{c+d-1-l}|\nabla b|^{2+\beta} e^{-f}  + \\
& & b^{c+d-l}\lp \nabla|\nabla b|^\beta, \nabla b\rp e^{-f} +b^{c+d-l}|\nabla b|^\beta (\Delta_f b) e^{-f},
\ean
plugging $\Delta_f|\nabla b|^\beta$ and $\Delta_f b$ with  (\ref{Delta-f-grad-b-alpha}) and (\ref{Delta-f-b}) to the above, we get
\ban
\lefteqn{ g(r_2)-g(r_1)} \\
&=&\int_{r_1\leq b\leq r_2}\bigg\{\frac{\beta}{4}b^{c+d-1-l}|\nabla b|^{\beta-2}\Big[\left|\Hess\, b^2\right|^2+\Ric_f(\nabla b^2,\nabla b^2)-4k|\nabla b|^4 +\nonumber\\
& &4(\beta-2)b^2|\nabla|\nabla b||^2\Big]e^{-f}+(2k-3+c+d-l)b^{c+d-l}\lp \nabla b, \nabla|\nabla b|^\beta \rp e^{-f} \bigg\}\nonumber\\
& &+(k+d-l)\int_{r_1\leq b\leq r_2}b^{c+d-l}\left[(c+d-l+k-1)b^{-1}|\nabla b|^{2+\beta}+\lp \nabla b, \nabla|\nabla b|^\beta \rp \right] e^{-f}.
\ean
This is  (\ref{A'-mono}) after grouping.
\qed

\sect{Monotonicity for $\Ric_f^N \ge 0$}
From Theorem~\ref{A-V'} and Theorem~\ref{rA'}, if $\Ric^N_f \ge 0$, we get many families of monotonicity quantities.

 Since $\Ric_f (\nabla b^2, \nabla b^2) = \Ric^N_f (\nabla b^2, \nabla b^2) + \frac{\lp \nabla b^2, \nabla f \rp^2}{N}$ and
 \be
|\Hess\,b^2|^2=\left|\Hess\,b^2-\frac{\Delta b^2}{n}g\right|^2+\frac{(\Delta b^2)^2}{n},   \label{Hess-b2}
\ee
we have
\ban
\left|\Hess\, b^2\right|^2 + \Ric_f (\nabla b^2, \nabla b^2) & \ge & \frac{(\Delta b^2)^2}{n} + \frac{\lp \nabla b^2, \nabla f \rp^2}{N} \\
& \ge & \frac{\left(\Delta_f b^2\right)^2}{n+N} = \frac{4k^2}{n+N} |\nabla b|^4.
\ean
Here we used the basic inequality $\frac{a^2}{p} + \frac{b^2}{q} \ge \frac{(a+b)^2}{p+q}$.

Therefore,  by Theorem~\ref{A-V'}, $(A^{\beta}_f -\alpha V^{\beta,p}_f)'(r) \ge 0$ if $\beta \ge 2$, $\lambda_1 \ge 0$ and $\lambda_2 +\frac{\beta k^2}{n+N} \ge 0.$ There are many solutions to these.

For example, if we let $\lambda_1=0$, namely $\alpha = 3k-p-l-2$, and $k \ge n+N$, then $\lambda_2 +\frac{\beta k^2}{n+N} \ge 0$ when $k\leq l\leq 2k-2$,which gives Theorem~\ref{A-V-mono}.

Letting $k=l=n+N$ in Theorem~\ref{A-V-mono}, we get
\begin{coro}  \label{A-V'k=l}
If $M^n (n\ge 3)$ has $\Ric^N_f \ge 0$, then, for $k=l=n+N$, $\beta \ge 2$, $p<n+N-\frac{\beta N}{n-2}$, and $0<r_1 < r_2$,
\[
(A^{\beta}_f - (2n+2N-p-2)  V^{\beta,p}_f)(r_2) \ge (A^{\beta}_f - (2n+2N-p-2)  V^{\beta,p}_f)(r_1).
\]
\end{coro}

Similarly, by Theorem~\ref{rA'}, $g(r) = r^c \left(r^d A_f^\beta (r) \right)'$ is nondecreasing if $\beta \ge 2, \lambda_3 +\frac{4 k^2}{n+N} \ge 0$ and $ \lambda_4 =0$. Again there are many solutions. $\lambda_4 =0$ requires that $c+2d = 3-3k+2l$. When $c=k-1, d=l-2k+2$ or $ c= 3-3k+2l, d=0$, and $k=l = n+N$, we have
\begin{prop} If $M^n (n\ge 3)$ has $\Ric^N_f \ge 0$, then for $0< r_1 < r_2$, $k=l= n+N$,
\ban
\lefteqn{ r_2^{k-1}(r^{2-k}A_f^\beta)'(r_2)-r_1^{k-1}(r^{2-k}A_f^\beta)'(r_1) }\nonumber\\
&\ge &\int_{r_1\leq b\leq r_2}\frac{\beta}{4}b^{-k}|\nabla b|^{\beta-2}\left\{\left|\Hess\,b^2-\frac{\Delta b^2}{n}g\right|^2 +
4(\beta-2)b^2|\nabla|\nabla b||^2\right\}e^{-f},
\ean
and
\ba
\lefteqn{ r_2^{3-k}(A_f^\beta)'(r_2)-r_1^{3-k}(A_f^\beta)'(r_1) }\nonumber\\
&\ge &\int_{r_1\leq b\leq r_2}\frac{\beta}{4}b^{2-2k}|\nabla b|^{\beta-2}\left\{\left|\Hess\,b^2-\frac{\Delta b^2}{n}g\right|^2 +
4(\beta-2)b^2|\nabla|\nabla b||^2\right\}e^{-f}. \nonumber \\
&& \label{rkA'k=l}
\ea
\end{prop}

Again, when $N=0$, $\beta =2$, these are the second and third  monotonicity formula in \cite{Colding}.

Now we prove Theorem~\ref{A-mono} which we restate it here.
\begin{theo}
If $M^n (n\ge 3)$ has $\Ric^N_f \ge 0$, then for $\beta \ge 2$, $k=l= n+N$, $(A^\beta_f)'(r) \le 0$ and $(V_f^{\beta, p})'(r) \le 0$ for $p<n+N-\frac{\beta N}{n-2}$. In fact \be
(A^\beta_f)'(r) \le - \frac{\beta}{4} r^{k-3} \int_{ b\ge r} b^{2-2k}|\nabla b|^{\beta-2} \left|\Hess\,b^2-\frac{\Delta b^2}{n}g\right|^2 e^{-f}.  \label{A'k=l}
\ee
\end{theo}
\Pf First we show $A_f^\beta(r)$ is bounded as $r \ra \infty$. By (\ref{nabla-b}), $|\nabla b(y)|\leq C(n, N, r_0)$ on $M \setminus B(x_0, r_0)$.
Hence for  $r\geq r_0$,
\ba
A^\beta_f(r)&=&r^{1-k}\int_{b=r}|\nabla b|^{1+\beta}e^{-f}\nonumber\\
            &\leq& C(n, N, r_0, \beta)r^{1-k}\int_{b=r}|\nabla b|e^{-f}. \nonumber
\ea
Define $h(r)=r^{1-k}\int_{b=r}|\nabla b|e^{-f}$, by (\ref{Iu'})
\ba
h'(r)=r^{1-k}\int_{b = r}\lp \nabla 1, \nu \rp e^{-f}=0. \nonumber
\ea
By (\ref{A-0}), $\lim_{r \ra 0} h(r) = \frac{n-2}{k-2} \vol(\partial B_1(0))e^{-f(x_0)}$.
Therefore
\be
A^\beta_f(r)\leq C(n, N, r_0, \beta)e^{-f(x_0)} \ \ \mbox{for} \   r\geq r_0.  \label{A-bounded}
\ee

Next we prove that $(A^\beta_f)'(r)\leq 0$. By (\ref{rkA'k=l}),  for $0<r_1<r_2$
\ba
r_2^{3-k}(A^\beta_f)'(r_2)\geq r_1^{3-k}(A^\beta_f)'(r_1). \nonumber
\ea
If there is some $\bar{r}>0$ such that $(A^\beta_f)'(\bar{r})>0$, then for all $r \ge \bar{r}$,
\ba
(A^\beta_f)'(r)\geq \left(\frac{r}{\overline{r}}\right)^{k-3}(A^\beta_f)'(\bar{r}) \ge (A^\beta_f)'(\bar{r}) >0. \nonumber
\ea
Namely $A^\beta_f(r)\rightarrow \infty$ as $r\rightarrow\infty$, contradicting to (\ref{A-bounded}).

To show $(V_f^{\beta, p})'(r) \le 0$, note that  by Corollary~\ref{A-V'-k=l}, $(A^{\beta}_f - (2n+2N-p-2)  V^{\beta,p}_f)'(r) \ge 0$. Hence $$(2n+2N-p-2)  (V^{\beta,p}_f)'(r) \le (A_f^\beta)'(r) \le0.$$

Now (\ref{A'k=l}) follows from (\ref{rkA'k=l}) since the fact that $(A^\beta_f)'(r) \le 0$ and $A_f^\beta(r)$ is bounded as $r \ra \infty$ imply there are sequence $r_j \ra \infty$ such that $r_j^{3-k}(A^\beta_f)'(r_j) \ra 0$.
\qed

As in \cite{Colding-Minicozzi}, one can decompose $\left|\Hess\,b^2-\frac{\Delta b^2}{n}g\right|^2$ further more to allow smaller $\beta$. Denote
\be
B = \Hess\,b^2-\frac{\Delta b^2}{n}g,  \label{B}  \ee
the trace free symmetric bilinear form, where $g$ is the Riemannian metric. One can decompose $B$ into the normal and tangential components. Let $B_0, g_0$ be the restriction of $B, g$ to the level set of $b$, define $B(\nu)$ the vector such that $\lp B(\nu), v\rp = B(\nu, v)$, and $B(\nu)^T$ its tangential component, so $|B(\nu)|^2 = (B(\nu, \nu))^2 + \left|B(\nu)^T\right|^2$. Then
\ba
|B|^2 &=&  (B(\nu, \nu))^2 + 2\left|B(\nu)^T\right|^2 +|B_0|^2 \nonumber \\
& = & |B(\nu)|^2 + |B(\nu)^T|^2 + \frac{(\tr B_0)^2}{n-1} + \left|B_0 - \frac{\tr B_0}{n-1}g_0\right|^2 \nonumber \\
& = & \frac{n}{n-1}|B(\nu)|^2 + \frac{n-2}{n-1} \left|B(\nu)^T\right|^2 +\left|B_0 - \frac{\tr B_0}{n-1}g_0\right|^2. \label{|B|^2}
\ea
Here we used the fact that $B$ is trace free so $\tr B_0 = - B(\nu, \nu)$.
Now
\begin{lem}
\be
\left|B(\nu)\right|^2 = 4b^2 |\nabla |\nabla b||^2 + \lambda^2 + 4\lambda b \lp \nabla |\nabla b|, \nu \rp,  \label{B-nu}
\ee
where $\lambda = 2|\nabla b|^2 - \frac{\Delta b^2}{n}= 2(1-\frac kn)|\nabla b|^2 - \frac 1n \lp \nabla b^2, \nabla f\rp$.
\end{lem}
\Pf Since \[
\Hess\, b^2 = 2b \Hess\, b + 2 \nabla b \otimes \nabla b, \]
we have
\ban
B(\nabla b) & = & 2b \Hess\, b (\nabla b) +2 |\nabla b|^2 \nabla b - \frac{\Delta b^2}{n} \nabla b \\
& = & b \nabla |\nabla b|^2  + \lambda \nabla b,
\ean
and
\[
B(\nu) = 2b \nabla |\nabla b| + \lambda \nu.\]
Hence \[
\left|B(\nu)\right|^2 = 4b^2 |\nabla |\nabla b||^2 + \lambda^2 + 4\lambda b \lp \nabla |\nabla b|, \nu \rp.\]
\qed

Therefore we have

\Rk  When $\lambda = 0$ (as in the case in \cite{Colding-Minicozzi}) or $\lp \nabla |\nabla b|, \nu \rp =0$, $|B|^2 \ge \frac{4n}{n-1}b^2 |\nabla |\nabla b||^2$. We also get all above monotonicity for $\beta \ge 1-\frac{1}{n-1}$ instead of $\beta\ge 2$.

\sect{Monotonicity for $\Ric_f \ge 0$}

When one only assumes  $\Ric_f \ge 0$, one does not have the extra term $\frac{\lp \nabla b^2, \nabla f \rp^2}{N}$ to combine with  $\left|\Hess\,b^2\right|^2$  to get $(\Delta_f b^2)^2$. First we study the monotonicity without using $\left|\Hess\,b^2\right|^2$.  Interestingly we still get many monotonic formulas.
\begin{theo} \label{A-V-mono-infty-gen}
If $M^n (n\ge 3)$ has $\Ric_f \ge 0$,  then for $\beta \ge 2$, $k,p$ such that $C(n,k,p) >0, (k-2)^2-4\beta k\ge 0$, $l_1\leq l\leq l_2$, where $l_1$ and $l_2$ are the solutions of $l^2+(2-3k)l+2k^2-2k+\beta k=0$, and $\alpha = 3k-p-l-2$,
we have $(A^{\beta}_f -\alpha V^{\beta,p}_f)'(r) \ge 0$.
\end{theo}
\Pf From (\ref{A-V'}), we only need to make sure $\lambda_1 \ge 0, \lambda_2 \ge 0$. By the choice of $\alpha$ we have $\lambda_1 = 0$ and $\lambda_2  \ge 0 $ if and only if $l^2+(2-3k)l+2k^2-2k+\beta k\leq 0$. This has real solution when $(k-2)^2-4\beta k \ge 0$.
\qed

 When $\beta =2, p=0, n \ge 4$, then for any $k\ge 12$, $l=\frac 32 k -1$ is between the $l_1,l_2$ occurred above, which is the first part of Theorem~\ref{A-V-mono-infty}.

 Similarly, by Theorem~\ref{rA'}, if $\beta = 2, c=d=0, l = \frac 32(k-1), k \ge 12$, then $\lambda_4 =0, \lambda_3 >0$, we get the second statement in  Theorem~\ref{A-V-mono-infty}.

To get monotonicity for $\Ric_f \ge 0$  when $k=l$, we need to use the $\left|\Hess\,b^2\right|^2$ term. Here are some formulas when $k=l=n$, which recover the formulas in \cite{Colding-Minicozzi} when $f$ is constant and $p=2$.
\begin{theo}  \label{mono-k-l-n}
When $k=l=n, p<n$, we have
\ba
\lefteqn{(A^{\beta}_f -(2n-p-2) V^{\beta,p}_f)'(r)  }  \nonumber\\
& & = r^{p-1-n} \int_{b\leq r} \frac{\beta |\nabla b|^{\beta -2}}{4b^{p}} \bigg\{ \left| B \right|^2 + \Ric_f (\nabla b^2, \nabla b^2)  \nonumber\\
& &  +4(\beta-2)b^2|\nabla|\nabla b||^2  + 4 |\nabla b|^2 \lp \nabla b^2, \nabla f \rp +  \frac{\lp \nabla b^2, \nabla f \rp^2}{n} \bigg\}e^{-f}, \label{A-V'-k=l=n}
\ea
where $B$ is given in (\ref{B}).
\end{theo}
\Pf From (\ref{A-V'-k=l})
\ban
\lefteqn{(A^{\beta}_f -(2n-p-2) V^{\beta,p}_f)'(r) = r^{p-1-n} \int_{b\leq r} \frac{\beta |\nabla b|^{\beta -2}}{4b^{p} }
\bigg\{ \left|\Hess\, b^2\right|^2  }   \nonumber\\
& & + \Ric_f (\nabla b^2, \nabla b^2) +4(\beta-2)b^2|\nabla|\nabla b||^2 - 4n |\nabla b|^4 \bigg\}e^{-f}.
\ean
By (\ref{Hess-b2}) and (\ref{Delta-f-b-2}),
\ban
\left|\Hess\, b^2\right|^2  & = & \left| B \right|^2+\frac{(\Delta_f b^2 + \lp \nabla b^2, \nabla f \rp)^2}{n} \\
& = & \left|B \right|^2+4n |\nabla b|^4 + 4 |\nabla b|^2 \lp \nabla b^2, \nabla f \rp +  \frac{\lp \nabla b^2, \nabla f \rp^2}{n}.
\ean
Plug this into above gives (\ref{A-V'-k=l=n}).
\qed

Similarly, letting $k=l=n, \ d=0,\ c = 3-n$ in Theorem~\ref{rA'} gives
\begin{theo}  \label{rA'-infty}  When $k=l=n$, for $0<r_1 < r_2$,
\ban
\lefteqn{ r_2^{3-n}(A_f^\beta)'(r_2)-r_1^{3-n}(A_f^\beta)'(r_1) }\nonumber\\
& = &\int_{r_1\leq b\leq r_2}\frac{\beta}{4}b^{2-2n}|\nabla b|^{\beta-2}\bigg\{\left| B \right|^2 +  \Ric_f (\nabla b^2, \nabla b^2) \nonumber \\
& & + 4(\beta-2)b^2|\nabla|\nabla b||^2 + 4 |\nabla b|^2 \lp \nabla b^2, \nabla f \rp +  \frac{\lp \nabla b^2, \nabla f \rp^2}{n} \bigg\}e^{-f}.
\ean
\end{theo}

\begin{coro} \label{mono-k=l=n}
If $\Ric_f \ge 0$, $\beta \ge 2$, $k=l=n, \ p<n$,  and \be
\left| B \right|^2 + 4 |\nabla b|^2 \lp \nabla b^2, \nabla f \rp +  \frac{\lp \nabla b^2, \nabla f \rp^2}{n} \ge 0,  \label{B-c}
\ee
then
\[
(A^{\beta}_f -(2n-p-2) V^{\beta,p}_f)'(r) \ge 0 \] and for  $0<r_1 < r_2$,
\[
r_2^{3-n}(A_f^\beta)'(r_2)  \ge r_1^{3-n}(A_f^\beta)'(r_1) .\]
\end{coro}

Note (\ref{B-c}) holds in particular if $\lp \nabla b^2, \nabla f \rp  \ge 0$ or $\lp \nabla b^2, \nabla f \rp  \le -4n |\nabla b|^2$. In general $\left| B \right|^2 + 4 |\nabla b|^2 \lp \nabla b^2, \nabla f \rp +  \frac{\lp \nabla b^2, \nabla f \rp^2}{n}$ could be negative and monotonicity in Corollary~\ref{mono-k=l=n} could even be reversed.  We illustrate this in the following example.
\begin{ex}  \label{Bryant-soliton}
Bryant soliton is a rotationally symmetric steady gradient Ricci soliton. It is $\mathbb R^n (n\ge 3)$ with the metric
\ba
g=dr^2+\phi(r)^2 g_{S^{n-1}}
\ea
where, as $r \ra \infty$,  $C^{-1}r^{1/2}\leq \phi(r)\leq C r^{1/2}$ and $\phi'(r)=O(r^{-1/2})$  for some positive constant $C$.
And the potential function $f(r)=-r + O(\ln r)$ as $r \ra \infty$. See e.g. \cite{Chow-etc}.

Then for $r > >0$,  $\left| B \right|^2 + 4 |\nabla b|^2 \lp \nabla b^2, \nabla f \rp +  \frac{\lp \nabla b^2, \nabla f \rp^2}{n}$  is positive when $n= 3$ and negative when $n \ge 5$, showing that  for  all $r_1 \le r_2$ big,
\ban
r_2^{3-n}(A_f^\beta)'(r_2)  & \ge &  r_1^{3-n}(A_f^\beta)'(r_1) \ \mbox{if} \ n=3,  \\
r_2^{3-n}(A_f^\beta)'(r_2)  & \le &  r_1^{3-n}(A_f^\beta)'(r_1) \ \mbox{if} \ n \ge 5.
\ean
\end{ex}
\Pf   The Green's function $G$ of the $f$-Laplacian with pole at the origin is $G=G(r)$, with
\[
G''+ (n-1)\frac{\phi'}{\phi} G'-f'G'=0.
\]
Hence $G' = -e^{f}\phi^{1-n}$, and $G(r)=\int_r^\infty e^{f}\phi^{1-n} ds$, both go to $0$ as $r \ra \infty$.

With $b=G^{\frac{1}{2-n}}$, we have $\frac{b' }{b} = \frac{1}{2-n} \frac{G'}{G}$, and
\ban \frac{b'' }{b'} &  = &   \frac{G''}{G'} - \frac{G'}{G} + \frac{b'}{b}  \\
& = &  f' -(n-1) \frac{\phi'}{\phi} -  \frac{(n-1)G'}{(n-2)G}.
\ean
Combine this with (\ref{|B|^2}) and (\ref{B-nu}),  we have
\ban
|B|^2&=& \frac{n}{n-1} |B(\nu)|^2 = \frac{n}{n-1} \left|2b \nabla |\nabla b| - \frac 1n \lp \nabla b^2, \nabla f\rp \nu\right|^2 \\
& = &  \frac{n}{n-1} \left( 2bb'' -\frac 1n 2bb'f' \right)^2 \\
& = & (2bb')^2 n(n-1) \left(\frac{G'}{(2-n)G} + \frac 1n f' - \frac{\phi'}{\phi} \right)^2 .
\ean
Now
\[
4|\nabla b|^2\lp\nabla b^2, \nabla f\rp = (2bb')^2 2f' \frac{b'}{b} = (2bb')^2 \frac{2}{2-n} f' \frac{G'}{G}.
\]
Hence
\ban
\lefteqn{ \left| B \right|^2 + 4 |\nabla b|^2 \lp \nabla b^2, \nabla f \rp +  \frac{\lp \nabla b^2, \nabla f \rp^2}{n} } \\
 & = &  (2bb')^2 \left[ n(n-1) \left(\frac{G'}{(2-n)G} + \frac 1n f' - \frac{\phi'}{\phi} \right)^2 - \frac{2}{n-2} f' \frac{G'}{G} + \frac{(f')^2}{n} \right].
 \ean
Since $\lim_{r \ra \infty} f'(r) = -1, \ \lim_{r \ra \infty} \frac{\phi'}{\phi}=0$ and
using L'Hospital's rule
\[
\lim_{r \ra \infty} \frac{G'}{G}=\lim_{r \ra \infty} \frac{G''}{G'}
               = \lim_{r \ra \infty} (f'+(1-n)\frac{\phi'}{\phi})=-1,  \nonumber
\]
we see that
\[ \lim_{r \ra \infty} \left[ n(n-1) \left(\frac{G'}{(2-n)G} + \frac 1n f' - \frac{\phi'}{\phi} \right)^2 - \frac{2}{n-2} f' \frac{G'}{G} + \frac{(f')^2}{n} \right]
=   \frac{-n^2+4n}{n(n-2)^2},\]
which is positive when $n =3$ and negative when $n\ge 5$.

Now the last claim follows from Theorem~\ref{rA'-infty}.
\qed

\end{document}